\documentclass[a4paper,10pt]{article}
\usepackage{hyperref}
\usepackage{amsfonts}
\usepackage{amssymb}
\usepackage[utf8]{inputenc}

\usepackage[english ]{babel}

\usepackage{bbding}

\usepackage[all]{xy}

\usepackage{tikz-cd} 
\usepackage[all]{xy}
\usepackage{amsfonts}
\usepackage{amssymb}
\usepackage{graphicx}
\usepackage{mathtools}
\usepackage{skak} 
\usepackage{foekfont}

\usepackage{calligra}
\usepackage[T1]{fontenc}

\usepackage{amsmath,mathdots}

\usepackage{stmaryrd}

\usepackage[utf8]{inputenc}

\usepackage{mathtools}
\usepackage{mathdots}
\usepackage{tikz}

\definecolor{amber(sae/ece)}{rgb}{1.0, 0.49, 0.0}
\newfont{\rsfsten}{rsfs10 scaled 1200}
\usepackage{MnSymbol}
\usepackage{tikz}
\usepackage{amscd}
\usepackage{MnSymbol}
\usepackage{wasysym}

\usepackage{mathrsfs}
\usepackage{pdfcolmk}

\usepackage{graphicx}
\newcommand*{\rom}[1]{\expandafter\@slowromancap\romannumeral #1@}

\usepackage{wrapfig}

\newcommand{\tightunderset}[2]{%
  \mathop{#2}\limits_{\vbox to .3ex{\kern-0.95ex\hbox{$#1$}\vss}}}
\newcommand{\tightoverset}[2]
{%
  \mathop{#2}\limits_{\vbox to .3ex{\kern-0.95ex\hbox{$#1$}\vss}}}

\newcommand{\oset}[2]{%
  {\mathop{#2}\limits^{\vbox to -.5\ex@{\kern-\tw@\ex@
   \hbox{\scriptsize #1}\vss}}}}
\usepackage[utf8]{inputenc}

\usepackage {ifsym}

\usepackage {url}
\usepackage{rotating}
\usepackage{mathrsfs}
\usepackage{amsbsy}

\title {  Isometric Immersions  with Controlled  Curvatures}
\author{Misha Gromov} %%%%%%%%%%%% 

\begin{document}
\maketitle

\begin{abstract} 
We  $\delta$-approximate strictly short (e.g. constant)  maps   between Riemannian manifolds $f_0:X^m\to Y^N$  for  $N\gg m^2/2$  by   $C^\infty$-smooth isometric immersions $f_\delta:X^m\to Y^N$ with curvatures 
$curv(f_\delta) < \frac {\sqrt 3}{\delta}$, for $\delta\to 0$.

  \end{abstract}
  \tableofcontents

  \section  { Compression and Approximation} 

%\subsection {Stable Isomeric Immersions} 
 A $C^\infty$-immersion of a smooth manifold $ X$ to a smooth Riemannian $Y=(Y, h) $,
$$f_0:   X\to Y$$
   is called  ${\cal II}_h$ or 
   just ${\cal II}$,   if the Riemannian metric inducing operator 
   from the space of smooth maps $f: X\to Y$ to the space $\mathcal G_+( X)$  of smooth semi definite positive quadratic forms on $X$, 
   $$\mathscr I= \mathscr I _h : \mathcal F=  C^\infty( X, Y)\to \mathcal G_+( X) \mbox { for }  f\overset {\mathscr I}\mapsto  g= f^\ast(h)$$    is {\it  \textbf infinitesimally \textbf  invertible in  a $C^\infty$-neighbourhood $\mathcal F_0\subset \mathcal F $   of $f_0$.}

This means that  the differential/linearization   of   ${\mathscr I}$,
   $$\mathscr L_f :T_f(\mathcal F) \to T_ {\mathscr I(f)} (\mathcal G), $$ 
   of $\mathscr I,\mbox { } f\in \mathcal F_0,$
    % $^\infty( X, Y)\to \mathcal G_+( X) \mbox { for }  f\overset {\mathscr I} L_f :(\phi} $
    is {\it   right invertible} by a {\it differential operator}
   $$\mathscr M_f:   T_ {\mathscr I(f)} (\mathcal G) \to \mathscr T_f(\mathcal F), \mbox  { }\mathscr    L_f\circ \mathscr  M_f=Id:  T_ {\mathscr I(f)}(\mathcal G)  \to  T_ {\mathscr I(f)}(\mathcal G),   $$
where $f\mapsto \mathscr M_f$ is a (possibly non-linear)  
differential  operator  defined on $\mathcal F_0$. 

If $Y=\mathbb R^N$, (and  in local coordinates for all $Y$, in general) the operator $\mathscr L$ 
can be written as an operator  on maps 
$\overrightarrow f:X\to Y$, 
 $$\mathscr L_{f}(\overrightarrow  f)= \mathscr I(f+\epsilon \overrightarrow  f)- 
  \mathscr I(f)+o(\epsilon), \hspace {1mm} \epsilon \to 0,$$ 
 and where 
 $\mathscr  M_{f}(\overrightarrow g)$ is a differential operator in  
 $(f, \overrightarrow g)$, which is linear in  $\overrightarrow g$ and which satisfies 
  $$\mathscr L_{f}(\mathscr  M_{f}(\overrightarrow g))= 
  \overrightarrow g.$$

{\it "Free" Example. }   
   Free  immersions $f$, i.e. where the second osculating spaces  $osc_2(f(x))\in T_{f(x)}(Y)$   have {\it maximal  possible}  dimensions, 
   $$dim(osc_2(f(x))= \frac {dim(X)(dim(X)+1)}{2}+dim(X)$$ at all points 
$x\in X$, are ${\cal II}$  
by the {\it Janet-Burstin-Nash 
Lemma.}

Consequently,  

\hspace {10mm}{\it generic $f$ are   ${\cal II}$ for $dim(Y) \geq 
   \frac {dim(X)(dim(X)+1)}{2}+ 2dim(X) $.}\footnote {See [Gr1986],  [Gr2017] and references therein.}

\textbf {1.A.  Definitions of  $m$-Free and of   Flat  $\mathcal{ II}^{[m]}$  Maps. } A  smooth immersion  $f:X=X^n\to Y$ is  $m$-free,   $m\leq n=dim(X)$,  if the  restrictions of $f$  to all $m$-dimensional  submanifolds in $X$ are free.

For instance, 
$$if  \hspace {2mm}N=dim(Y)\geq {m(m+1)\over 2} +m(n-m)+2n,$$
{\it then  generic $f$ are $m$-free.}\vspace{1mm}

  An  immersion  $f:X=X^n\to Y$ is   {\it flat 
$\mathcal{ II}^{[m]}$},   
 $m\leq  n$, if the induced metric in $ X$ is Riemannian flat, i.e. locally isometric to $\mathbb R^n$, and if the restrictions  of $f$  to all  {\it flat}, (i.e. locally isometric to $\mathbb R^m\subset \mathbb R^{ n} $) $m$-dimensional   submanifolds  $X^m\subset Y$  are   $\mathcal{ II}$. 

For instance, isometric $m$-free immersions 
 of  flat tori  are    flat 
$\mathcal{ II}^{[m]}$.

Such immersions, especially of the split tori $\mathbb T^n=\underset {n}{\underbrace{\mathbb T^1\times...\times \mathbb T ^1}}$ to the Euclidean spaces, play a  special  role in our arguments,\footnote { This similar to how it is with non-isometric immersions with controlled curvature  
 in [Gr2022].}  where we use below the following.

\textbf {1.B. Example}. {\sf  Let $ \mathbb T^n$ be the torus with a flat Riemannian  metric, i.e. the universal covering of  this  $\mathbb T^n$ is isometric to
$\mathbb R^n.$} Then:

{\sl $\mathbb T^n$ admits a   free isometric $C^\infty$-immersion
to $\mathbb R^{\frac{n(n+1)}{2}+n+2}$}

 and 

{\sl an  isometric  $\cal II$-immersion to
$\mathbb R^{\frac{n(n+1)}{2}+n+1}.$}\vspace {1mm}

\textbf {1.C.} {\it Remarks on the Proof.} (a)  The existence of free isometric immersions   
of flat  split  $n$-tori to $N$-dimensional Riemannian manifolds 
is proven for $N\geq \frac{n(n+1)}{2}+n+2$ in section 3.1.8 in [Gr1986] 
and since non-split flat tori can be approximated by  finite coverings of split ones, 
the general case follows from the Nash implicit function theorem.

(b)  The  existence  of  $\cal II$-immersion   of split  tori to $\mathbb R^{\frac{n(n+1)}{2}+n+1}$
  is (implicitly) indicated in exercise  on  p. 251 in  [Gr1986]\footnote {We explain this in section 3. Also notice  that there is no known obstructions   to the existence of free isometric immersions  of
  flat  $n$-tori    to $\mathbb R^{\frac{n(n+1)}{2}+n}$, but no example  of a free (isometric or not) immersion from $\mathbb T^n $ to   $\mathbb R^{\frac{n(n+1)}{2}+n}$  for $n\geq 2$  had been found  either.}, while the recent result  by DeLeo [DeL2019] points toward a  similar possibility   for 
$$dim(Y)\geq   \frac{n(n+1)}{2}+n-\sqrt \frac {n}{2}+\frac {1}{2}. $$

In fact, it is not impossible, that such immersions exist for 
$$dim(Y)\geq\frac {n(n+1)}{2}+1, $$ 
and it seems not hard to show these don't exist for $dim(Y)\leq   \frac {n(n+1)}{2}$.

\vspace {1mm}
    
     \textbf {Definition of $\bf curv(f)= curv(f(X))$.} This is (as in [Gr 2022]) the curvature of a manifold     $X$ immersed by $f$  to a Riemannian  $Y$, that   is
  {\it the supremum of the "$Y$-curvatures" }, i.e. curvatures measured in the Riemannin geometry of $Y\supset X$,    of geodesics $\gamma\subset X$, 
    for the induced Riemannin metric in $X$,  
   $$ curv(f)=curv(f(X))=\sup_{\gamma\in X} curv_Y(\gamma).$$

\textbf {1.D.   $\sqrt 3$-Remark.}   One has only a limited control  over the curvatures of the above immersions $\mathbb T^n\to Y$, even   for $Y=\mathbb R^N$, but 
 we shall prove  in the next section  the existence of  
 
 {\it a free  isometric   immersion from  the $n$-torus  to the unit $N$-ball for $N\geq  \frac{n(n+1)}{2}+2n$   with the curvature bounded by a constant $D$ independent of $n$,} where, conceivably, $D< \sqrt 3.$

In fact, $ \sqrt 3$ is {\it asymptotically optimal}, since, according to {\it  Petrunin's inequality},\footnote { See \url {https://anton-petrunin.github.io/twist/twisting.pdf} and [P 2023].} \vspace {1mm}

{\it smooth isometric immersions $f:\mathbb T^n\to B^N(1)$ satisfy
 for all $n$ and $N$: 
$$curv(f)\geq \sqrt {3{n\over  n+2}}.$$}

\textbf {1.E.   Compression Lemma.} 
{\sf Let $$ F : \mathbb T^{ n}\hookrightarrow B^{ N}(1)\subset \mathbb R^{ N}$$  be a
 flat 
$\mathcal{ II}^{[m]}$-immersion.

Let $X^m=(X^m,g)$ be a compact   Riemannian $m$-manifold, possibly with a  boundary, which admits a smooth $(1+\varepsilon)$-bi-Lipschitz immersion  
$\phi_\varepsilon:  X^m\to  \mathbb T^{n}$.
If $\varepsilon  \leq  \varepsilon_0( m)>0$. \footnote {Possibly, this
$\varepsilon_0$ doesn't depend on $m$.}

Then
 there exist $C^\infty$-smooth isometric immersions 
$$f^\circ_i: X\hookrightarrow B^{ N}\left (\frac {1}{i} \right ), i=1,2,..., $$  
such 
that 
$$curv_{f^\circ_i}(X)\leq i\cdot curv(F(\mathbb T^{n}))+O(1).$$}
 
{\it Proof.}  Compose the maps $F$ and $\phi_\varepsilon $ with the homothetic endomorphism $t\mapsto i\cdot t$ of the torus, 
$$X  \overset {\phi_\varepsilon }\to \mathbb T^{n} \overset {i\cdot}   \to \mathbb T^{n}\overset {F}\to B^N(1) $$
 and observe that the resulting composed  maps, say 
 $$f_{i,\varepsilon} :X\to B^N(1)$$  
 satisfy the following conditions.

$\bullet_{\varepsilon} $  The map $f_{i,\varepsilon}$ is a  smooth $(1+\varepsilon)$-bi-Lipschitz immersion with respect the metrics 
$i^2g$.

$\bullet_{1/i} $ The covariant derivatives of the metric $i^2g$ and the covariant derivatives of the 
 induced metric $f^\ast (h_{Eucl})$ with respect to $i^2g$ converge to zero for $i\to \infty$, 
 $$\max (||\nabla_{ig}^j (i^2g||,  ||\nabla_{ig}^j (f_{i,\varepsilon} ^\ast (h_{Eucl})||)\leq const_ji^{-j}. $$

$\bullet_{\cal II} $ The immersions $f_{i,\varepsilon}$ are {\it uniformly $\mathcal{II}_{h_{Eucl}}$}:\vspace {1mm}

  {\sf the $j^2g$-covariant derivatives of the
  (the coefficients of the) differential operators $\mathscr M= \mathscr M_i=\mathscr M_{f_{i,\varepsilon}}$,
  which
  invert the linearized operator $\mathscr I: f_{i,\varepsilon}\to \mapsto   f^\ast_{i,\varepsilon}(h_{Eucl})$ on $X$ are  bounded,
for all $i$ independently of $i$,
$$||\nabla^j_{g_{i,\varepsilon}}(\mathscr M)||\leq const=const_{m,j}.$$}

It follows from the (generalized) Nash implicit function theorem  
 (section 2.7.2 in[Gr1986]) that, for  a sufficiently small $\varepsilon>0$, depending only on $m$, the maps $f_{i,\varepsilon} $ for sufficiently large $i$ admit $C^\infty$-small, convergent to $0$ with $\varepsilon \to 0$,  perturbations   to smooth isometric immersions   
 $ f_i^{\Circle}: (X,i^2g)\to B^N(1)$.
 
 Since the curvatures of these immersions  are bounded by the    curvature of $F$ plus $O(1/i)$, the immersions 
 $$f^{\circ}_i=i^{-1}  f_i^{\Circle}: (X,i^2g)\to B^N(1/i)$$
 are the required ones. QED.
 
  \textbf{1.F.  Local Compression Corollary.} {\sf Let $X^m$ be a smooth Riemannian manifold with the sectional curvature bounded by 
  $$|sect.curv(X^m)|\leq 1$$ and where  the injectivity radius at a given point is bounded from below
  $$inj.rad_ {x_0}(X^m)\geq 1.$$
 Then there exists a  constant $\rho=\rho_m>0$,\footnote {Probably, what we say here holds for
 $\rho\geq 10^{-10}$  and all $m$.} such that the ball $B(\rho)=B_{x_0}(\rho)\subset X$
 admits a {\it smooth isometric immersion} to the Euclidean space
 $$f:B(\rho)\to \mathbb R^{\frac {m(m+1)} {2} +m+1}.$$
 Moreover there exist immersions 
 $f_\delta: B(\rho)\to \mathbb R^{\frac {m(m_1} {2} +m+1}$
 for all $0<\delta\leq 1$  with the {\it diameters of the images bounded by} 
 $$diam (f_\delta(B(\rho))\leq \delta$$
   and the{\it  curvatures of these images bounded by}  
 $$curv( f_\delta(B_{x_0}(\rho))\leq \frac {C_m}{\delta}.\footnote{The proof of 1.B in [Gr1986]  for $Y=\mathbb R^{\frac {m(m+1)} {2} +m+2}$ shows that $C_m< (100m)^{100m}$.} $$}

{\it Proof.} The assumptions on the curvature  and  the injectivity radius of $X$ imply that the $B_{x_0}(\rho))$ is 
$(1+3\rho)$-biLipschitz to the $\rho$-ball in the flat torus.

{\it Remark}.  We shall proof  in the next section the existence of 
 smooth isometric  $f_\delta: B_{x_0}(\rho)\to \mathbb  R^{\frac {m(m+1)} {2} +2m+1}$ with $diam (f\delta_{x_0})\leq \delta$  and the curvatures of these images bounded by  
 $$curv( f_\delta(B_{x_0}(\rho))\leq \frac {C}{\delta}$$
 for a universal constant $C$
 \footnote{Probbaly $C<100$.}

 %{\sf Let $$ f : \mathbb T^{ n}\hookrightarrow B^{ N}(1)\subset \mathbb R^{ N}$$  be a
 %flat 
%$\mathcal{ II}^{[m]}$-immersion.

 \textbf {1.G.   Approximation Lemma.} {\sf 
Let 
$$X^m=(X^m, g_\varepsilon)\overset { \phi_\varepsilon}\to \mathbb T^n \overset{ F} \to \mathbb R^N$$ 
be as in 1.E, let  $Y=(Y,h) $ be a smooth Riemannian manifold and let $f_0:X^m\to Y$  be a $C^\infty$-smooth 
map. 

Let the pullback to $X$ of  tangent bundle of $Y$,
$$f_0^\ast(T(Y))\to X,$$ 
admit $N${\it  independent vector fields normal to the image of the differential of $f_0$.}

For instance,  $dim(Y)\geq  N$  and $f_0$ is a constant map, or $f_0$ is an immersion homotopic to a constant map and 
$dim(Y)\geq N+2m-1$. 

If   $0<\varepsilon\leq \varepsilon_0  (m)$,
then, for all $i=1,2,...$, there exist a $\delta_i$-{\it approximation} of $f_0$ for $\delta_i \leq\frac {1}{i}$ by  $C^\infty$-immersions $f_i:X^m\to Y$   with
$$ curv_{f_i}(X)\leq i\cdot curv(F(\mathbb T^{ n})) +o\left({i}\right)$$
and where $f_i$  increase  the induced Riemannian metric $g_0=f^\ast_0(h)$ in $X$ by the above $g_\varepsilon$: 
$$f_i^\ast(h)=g_0+g_\varepsilon. \footnote {$
\delta_i$-{\it Approximation} signifies that $dist_Y(f_0(x), f_i(x))\leq \delta_i$,  $x\in X$.}.$$}

{\it Proof.} Let  $E=E_{N,\delta_0}:X\times  B^N(\delta)\to Y$ be the exponential map defined by the $N$ vector fields, where this map is defined for all  $\delta_0$ if $Y$ is complete $Y$ and if $Y$ is non-compete, then $E_{N,\delta_0}$ is defined if   the $\delta_0$-neighbourhood of  $f_0(X)\subset Y$ is compact.\footnote {We assume here that $Y$ has no boundary.}

Let $\delta_i\leq \delta_0$ and let us  restrict the map $E$ to the
 graph of the above map 
 $f^\circ_i:X\to B^N(\delta_i)\subset \mathbb R^N$,
 $$\Gamma_{i}=\Gamma_{f^\circ_i}:X\hookrightarrow  X\times B^N(\delta_i),$$
where our  $f^\circ _i$   is now isometric for the metric $g_\varepsilon$ on $X$.

Since the $N$ fields are {\it normal} to $f_0(X)\subset Y$,
the Riemannin metric in $X$  induced by $\Gamma_i\circ E: X\to Y$ is 
$(1+\varepsilon +const \cdot\delta_i)$-bi-Lipschitz to  $g_0+g_\varepsilon $;  hence,   the (generalized) Nash implicit function theorem applied as in the proof of 1.E, now   to the maps
 $\Gamma_i\circ E$ for small $\varepsilon>0$ and   large $i$,    delivers 
$C^\infty$-perturbations to these maps to the required immersions
$f_i:X\to Y$  with $f^\ast_i(h)=g_0+g_\varepsilon$
and with $ curv_{f_i}(X)\leq i\cdot curv_{ f}(\mathbb T^{ n}) +o\left({i}\right).$

\vspace {1mm}

 \textbf {1.H.  Global Approximation  Corollary.} {\sf Let $X^m=(X^m,g)$ and $Y^N=(Y^N, h)$  be  smooth  Riemannian manifolds and $f_0:X \to Y=(Y,h)$ be a smooth {\it strictly short  map}, i.e. the quadratic differential form $g-f^\ast(h)$ is positive definite. 
 
 If $X^m$ is compact, if $X$ admits a smooth immersion to  $\mathbb R^n$ and if
 $$N\geq{n(n+1)\over 2} +n+1, $$
 then there exists  $\delta_i$-{\it approximation} of $f_0$ for $\delta_i \leq\frac {1}{i}$, $i=1,2,...$, by  isometric $C^\infty$-immersions $f_i:X^m\to Y$   with
$$ curv(f_i((X))\leq i \cdot C_m+o\left({i}\right).$$}
 
 {\it Proof.} Apply the lemma to   smooth $(1+\varepsilon)$-bi-Lipschitz  immersions $\phi_\varepsilon: X^m\to\mathbb  T^n=\mathbb R^n/\mathbb Z^n$,  delivered  by the Nash-Kuiper $C^1$-immersion theorem and to an isometric  $\cal II$-immersion 
$ F:\mathbb T^n\to \mathbb R^{N=\frac{n(n+1)}{2}+n+1}$  from 1.B.
 
 {\it Remarks.} (a) Since all  $X^m$, $m\geq 2$,  admit  smooth immersions to $\mathbb R^{2m-1}$
 by the Whitney theorem, the inequality 
 $$N\geq {2m(2m-1)\over2} +2m\approx 2m^2 $$
 suffices for all $X^m$. 
 This is  the worst case  dimension-wise.

Our bound on $N$  is  much better  for  $n=m+1$, e.g. for compact hypersurfaces $X^m\subset \mathbb R^{m+1}$, where one needs
 $$N\geq {(m+1)(m+2)\over 2} +m+2={m(m+1)\over 2}+ 2m+3;$$
 but this still  seems far from  optimal. 
 
 The best we  can  get for smaller $N$, namely for 
 $$N\geq     \frac {m(m+1)}{2} +m+ 1,$$
 is the following  special result.

\textbf {1.I. Flat Torus Approximation Theorem.}  
 {\sf Let $\mathbb T^m_g$  the torus with an {\it invariant} (hence {\it flat}) Riemannian metric  $g$  and  let 
  $$f_0:\mathbb T_g^m\to Y$$
 be
 a smooth map, where $Y=(Y,h)$ is a 
 compact   $N$-dimensional $C^\infty$-smooth  Riemannian manifold, possibly with a boundary, such  that 
  $$dist (f_0(\mathbb T_g^m ), \partial Y)>0.$$

If $N=dim (Y) \geq \frac {m(m+1)}{2}+m+1$,
and if the  bundle  induced by   $f_0$ from the tangent bundle of $Y$, that is $f_0^\ast (T(Y)\to \mathbb T^m_g$ is  {\it trivial}
(e.g.  $f_0$ is contractible or  $Y$ is stably parallelizable), 
then, for all $\delta>0$,\footnote {Since we are concerned with $\delta \to 0$,  we assume here and below  that  $\delta\leq 1$.}  the map  $f_0$ admits a $\delta$-approximation by 
{\it metrically homothetic} maps  
  $$f_\delta:\mathbb T_g^m\to Y,$$
i.e. such that $$f^\ast(h)=\lambda \cdot g \mbox {  for some }\lambda>0$$
and where 
 $$curv(f_\delta(\mathbb T^m))  \leq \frac 
{const_m}{\delta} +o\left (\frac {1}{\delta}\right),$$ }

{\it Proof. }  Let 
  $E=E_{N,\delta_0}:X\times  B^N(\delta)\to Y$ be the exponential map, similar to that in the proof of 1.F but now {\it not required to be 
    normal} to $f_0(X)$.
  
 Let $F:\mathbb T^m\to B^N(1)$, $N=\frac {m(m+1)}{2} +m+1$ be as in 1.B
 and let  
 $$F_{ij} : \mathbb T^m\to B^N(1/i)\mbox  {  for }  t\mapsto \frac {1}{i} F(jt).  $$

  Let $\delta_i\leq \delta_0$ and let us  restrict the map $E$ to the
 graph of the map $F_{ij} $
 $$\Gamma_{F_{ij}}= \Gamma_{i,j}:X\hookrightarrow X\times B^N(\delta_i).$$

 Let  the ratio $\lambda = j/i$  be very  large.
Then,   in terms of the metric  $\lambda g$, the metric induced by $E\circ \Gamma_{ij}:\mathbb T^m\to Y$ becomes $C^\infty$-close to $\lambda g$ and 
the proof follows as in  1.E and 1.G  by the (generalized) Nash implicit. function theorem.

{\it Remarks.} The proof of   theorem(c)  on p. 294 in [Gr1986]   delivers  (a stronger version of) the  
"local compression"  for surfaces,   $X^2\to B^4(\delta)$, and this, seems  to imply   
 the torus approximation  theorem for  $\mathbb T^2\to Y^N$  and all $N\geq 4$.

We don't know  if one could    comparably   improve  bounds on $N$ in general,  but
 in the next section we prove an approximation theorem for all $X$ with  the constant $C$  {\it independent of $m$} and  with 
a  slightly  improved  dimension   bounds  in some cases.

{\it \textbf{Acknowledgments.}} I want to thank an anonymous referee  for  useful suggestions  and for pointing out  several errors in an  earlier version of the paper.

\section {Free  Isometric Imbeddings of Tori to the Unit Balls with Small Curvatures }

 Let
$$\mathbb T_{ \sf Cl}^N\subset S^{2N-1}\subset B^{2N}\subset \mathbb R^{2N}$$ 
be the Clifford torus, which observe, has  
the Euclidean curvature 
$$curv(\mathbb T_{\sf  Cl}^N\subset \mathbb R^{2N})=\sqrt N.$$

\textbf{ 2.A.  $\Delta(n,N)$-inequality.} If  
$1\leq  \frac {n^2}{2}\leq N$,  then there exits  a (flat invariant)
  subtorus $\mathbb T_0^n\subset \mathbb T_{ \sf Cl}^N$, 
 the Euclidean curvature of which satisfies
$$ curv (\mathbb T_0^n \subset \mathbb R^{2N})
 \leq  \sqrt 3 \sqrt \frac{n}{n+2}+\Delta(n,N)\footnote{The summand   $\sqrt 3 \sqrt \frac{n}{n+2}$ is optimal by Petrunin's inequality.}$$
where $\Delta(n,N)$ is bounded by a universal
constant, which, in fact, {\it vanishes}  for  $N\geq 8(n^2+n)$.\footnote { 
 Probably,  $\Delta(n,N)< 10$  for all $N\geq n(n-1)/2$  and, possibly  $\Delta(n,N)\leq 1/n$ for $N\geq \frac {n(n+2)}{2}$.}   

This follows from the $D(m,N)$-inequalities  2.1.E and 3.B in [Gr2022].

\textbf {2.B. $m$-Freedom  Corollary.} {\sf If $m\leq n$ and
$$ \frac {m(m+1)}{2} + m(n-m)+2n\leq 2N,$$
then, for all $ \varepsilon >0$, the   $n$-torus $ \mathbb T^n  $   admits an  immersion  $ \mathbb T^n\overset {F_\varepsilon} \hookrightarrow  B^{2N}\subset \mathbb R^{2N}$ such that   

$\bullet_{flat} $  the induced  Riemannian metric $F_\varepsilon^\ast(h_{Eucl})$ in  $ \mathbb T^n  $ is {\it Riemannin flat};

$\bullet_{curv}$ the Euclidean {\it curvature} of this  torus is {\it bounded} by $$curv(F_\varepsilon(\mathbb T^n)) \leq  \sqrt \frac{3n}{n+2}+\Delta(n,N)+\varepsilon;$$

$\bullet_{free}$ The restrictions of $F_\varepsilon$ to all $m$-dimensional submanifolds in $\mathbb T^n$   are {\it free}.} 

{\it Proof.} The required $F_\varepsilon$ is obtained by generically $\varepsilon$-perturbing the lift of the  above 
 $T_0^n\subset \mathbb T_{ \sf Cl}^N$ to a finite covering of $\mathbb T_{ \sf Cl}^N$  as follows.

Firstly, replace   $ \mathbb T_0^n$ by a generic (flat invariant) subtorus,   
$\mathbb T_\varepsilon^n\subset \mathbb T_{ \sf Cl}^N$ {\it tangentially} 
$\varepsilon$-close to $ \mathbb T_0^n$, i.e. where 
 the tangent spaces to $\mathbb T_\varepsilon^n$  are $\varepsilon$-parallel in $ \mathbb T_{ \sf Cl}^N$  to these of $T_0^n$,\footnote{This $ \mathbb  T_\varepsilon^n$ is far from $T_0^n$ as a subset in $\mathbb T_{ \sf Cl}^N$.} where such an 
 $ \mathbb  T_\varepsilon^n$ can be chosen split, i.e. being Riemannian product of $1$-tori.
 
 Moreover, we let  
 $$\tilde{ \mathbb T}_{ \sf Cl}^N=\underset {N} {\underbrace{\mathbb T_1^1\times ...\times  \mathbb  T^1_i \times ...\times \mathbb  T_N^1}},$$
 (the circles $ \mathbb  T^1_i$ may have different lengths $l_i$),
 be a  split  finite covering of the Clifford torus
 such that 
 $$\tilde {\mathbb T}^n=
 \underset {N} {\underbrace{\mathbb T_1^1\times ...\times  \mathbb  T^1_n}}$$
equals a covering of $\mathbb T_\varepsilon^n$.
 
 %\footnote {Split subtori in  $ \mathbb T_{ \sf Cl}^N$  are tangentially dense among all subtori.}

 Secondly,  
  %$$\mathbb  T_\varepsilon^n \times \mathbb R^{N-n} \to \mathbb  T_{ \sf Cl}^N$$ 
  %be the natural covering map and 
    perturb the embedding $\tilde {\mathbb  T}^n\subset \tilde { \mathbb  T}^N_{\sf Cl} $
keeping the induced metric Riemannin flat in $N-n$ steps, where, at each step, 
we perturb a split $K$-torus, which contains  $\tilde {\mathbb  T}^n$, in a $(K+1)$-torus, 
$$\tilde {\mathbb  T}^n\subset \mathbb  T^{K-1}\times \mathbb T_i^1\subset \mathbb  
T^{K-1}\times (\mathbb T_{i}^1\times \mathbb T^1_{i+1})=\mathbb T^{K+1}\subset 
 \tilde { \mathbb  T}^N_{\sf Cl}, \mbox { } K=n,...,N-1,  i=n+1,... N,$$
 by  approximating the embedding 
$\mathbb T_{i}^1 \subset \mathbb T_{i}^1\times \mathbb T^1_{i+1} $ by
a generic isometric embedding 
$(1+\epsilon)\cdot \mathbb T_i^1\overset {I_\varepsilon}\to  \mathbb T_i\times
 \mathbb T^1_{i+1}$.

Here $(1+\epsilon)\cdot\mathbb T_i^1$
is
the same  circle as  $\mathbb T_i$ but with the metric of total
 length equal to $(1+\epsilon)\cdot length ( \mathbb T_i)$ and  where, finally, we  let 
$$\mathbb  T^{K-1}\times \mathbb T_i^1\ni (\theta, t)\to (\theta,  I_\epsilon(t))\in \mathbb  T^{K-1}
\times \mathbb T_i^1\times \mathbb T_{i+1}^1.$$
 
If $ \frac {m(m+1)}{2} + m(n-m)+2n\leq 2N,$ and granted all was done "sufficiently generically", then the resulting map  
    $\tilde {\mathbb  T}^n \to B^{2N}(1)$
via  $\tilde { \mathbb  T}^N_{\sf Cl}\subset  B^{2N}(1)$
is free on all $X^m\subset \tilde {\mathbb  T}^n. $.

Checking this, which is similar to "{\it Making Non-free Maps Free}" on p. 259 in [Gr 1986], is left to the reader.

\textbf {2.C. Corollary.} 
{\sl  Let $X^m$ and $Y^M$  be  smooth  Riemannian manifold and $f_0:X \to Y$ 
 be a smooth strictly short  map.
 If $X^m$ is a compact manifold, which admits a smooth immersion to $\mathbb R^n$, and if
 $$M>  \frac {m(m+1)}{2} + m(n-m)+2n,$$
then $f_0$ can be     $\delta_i$-{\it approximated}, by  isometric $C^\infty$-immersions $f_i:X^m\to Y^M$,  for $ \delta_i \leq\frac {1}{i}$, $i=1,2,...$,   such that
$$  curv (f_i(X))
\leq i\left (\sqrt \frac{3n}{n+2}+\Delta(n, \lfloor {M/2}\rfloor)\right)+
o(i).$$}
 
 For instance, 
 
 {\it  if $X$ is a Euclidean hypersurface, then such an approximation is possible for $$M>  \frac {m(m+1)}{2} + 3m+2$$
 with 
 $$  curv (f_i(X))
\leq i\left (\sqrt \frac{3(m+1)}{m+3}+\Delta(m+1,\lfloor M/2\rfloor)\right )+o\left ({i}\right).$$
}

And -- this is the worst case -- the inequality  $$M> \frac {m(m+1)}{2}  + m(m-1)+
4m-2 $$ is sufficient for all  compact $ X^m$, where the maps $f_i$  satisfy: 
 
 $$  curv (f_i(X))
\leq i\left (\sqrt \frac{3(2m-1)}{2m+1}+\Delta(2m-1,\lfloor M/2 \rfloor)\right )+o\left ({i}\right).$$
\footnote{If the right-hand sides  of the  inequalities  $M>...$ are even, then these  may be  replaced   by $M\geq ....$ . }

\textbf  {2.D.  Immersions with  Prescribed Curvatures.}  
 The  {\it symmetric   normal curvature} $\Psi_f$  of an immersion $f:X\to Y$ is the   "square" of the second fundamental form, that is the
 symmetric differential 4-form, such that 
$$\Psi_f(\partial,\partial,\partial,\partial)=|| \nabla_{\partial,\partial}f||^2_Y$$
for all tangent vectors  $   \partial \in T(X).$

Observe that if $f$ is free, then  $\Psi_f$ is {\it  positive  definite}. This means that 

{\sf  the 4-form  $\Psi_x$ on the tangent space  $T=T_x(X)(=\mathbb R^m)$  is contained,  for all  $x\in X$, in the interior of the convex hull of the  $GL(m)$-orbit of the fourth power of a non-zero 1-form on $T$.\footnote {This interior makes the unique  open $GL(m)$-invariant (non-empty!) minimal convex cone  in the space 
$T^{\circledS 4}$ (of dimension $m(m+1) (m+2)(m+3)/ 24$, $m=dim (T)$  of  symmetric $4$-linear forms (4d-polynomials)   on $T$.

(This cone  is strictly smaller  than the cone of the forms $\Phi$, which  are positive  as polynomials, $\Phi (t,t,t,t )>0$, $t\neq 0$.)}  }

Also observe that
 $$\sup_{||\partial||=1}\Psi_f(\partial,\partial,\partial,\partial)=(curv(f))^2.$$

\textbf {2.C.  $C^2$-Curvature Theorem.}\footnote{See 3.1.5.(A) in [Gr 1986].}   {\it Let $X=X^m$ and $Y=Y^M$ be smooth Riemannian manifolds, let   $f:X \to Y$  a free isometric $C^\infty$-immersion and let $\Psi$ be a  {\it symmetric positive definite} differential 4-form on $X$.

 If 
$$M=dim(Y)\geq{m(m+1)\over 2} +3m+5,$$
then   $f$ 
can be  arbitrarily finely  $C^1$-approximated by isometric 
$C^2$-immersions $f'$  with the  increase of their normal curvatures by   $\Psi$.\footnote{In general, such an $f'$ can't be $C^\infty$   for large $m$ , but   possibly, these $C^2$-smooth $f'$ exist for all $m\geq 2$  and all $M$.} 
$$\langle ||\nabla_{\partial \partial} f'||^2=||\nabla_{\partial \partial}  f||^2 +\Psi(\partial,\partial,\partial,\partial)$$
for all tangent vectors $\partial \in T(X)$.}

 {\it \textbf {2.E. Euclidean Example.}}  {\it The standard embedding $f_0:\mathbb R^n\to\mathbb R^{(n + 2)(n + 5)/2}$ 
can be $C^1$-approximated by isometric $C^2$-embeddings $f$  with an arbitrary strictly positive definite normal curvature form 
$\Psi$.}

{\it Proof. }  $C^\infty$-approximate    $f_0$ by  
{\it free} isometric embeddings (as in 1.B) and apply 2.C.

{\it \textbf {2.F. Toric  Example.}} {\it Let  $\Psi$ be a smooth  symmetric positive definite differential  4-form on the $m$ torus $\mathbb T^m$. %and $ 3_{\mathbb T^m}$ be the constant form, such that $$3_{\mathbb T^m}(\partial,\partial,\partial,\partial)=3||\partial||^4$$

 If $M\geq 16(m^2+m) $, 
then  there exist a $C^2$-immersion $f:\mathbb T^M\to B^m(1)$, such that the induced metric is flat (split  if you wish) and 
$$\Psi_f(\partial,\partial,\partial,\partial)=\Psi(\partial,\partial,\partial,\partial)+\frac {3m}{m+2}||\partial||^4.$$ }

{\it Proof.} Observe that, according to 3.1 from [Gr 2022] (this also follows  from Petrunin's inequality),  the isomeric immersion  $f_0:\mathbb T^m\to B^M$
with $(curv(f_0))^2\leq {3m\over m+2}$ has constant curvature, 
$$\Psi_{f_0}((\partial,\partial,\partial,\partial)= {3m\over m+2}||\partial||^4,$$
and apply 2.B and  2.D.

\section {Perspective and Problems}

The $\mathcal {II}$-property of free immersions  was already  implicitly present 
in the  proof of the algebraic Janet  lemma (1926), where  this lemma brings the  
isometric immersion equations to the Cauchy–Kovalevskaya  form and thus implies that real analytic Riemannin $n$-manifolds locally $C^{an}$-immerse to $\mathbb  R^\frac {n(n+1)}{2}$.\footnote{The Cauchy–Kovalevskaya theorem also yields  a weak form of the  Nash real analytic implicit function theorem.  This,  by Janet's lemma combined  with an analytic version of  Nash twisting    argument,  delivers  isometric $C^{an}$-immersions of compact $C^{an}$-manifolds to Euclidean spaces, see [G-J 1971] and p 54 in Appendix 11 in  [G-R   1970].}
  \footnote{The existence of local isometric $C^\infty$-immersions of   $C^\infty$-manifolds $X^n\to \mathbb R^N $ is known for 
  $N= \frac {n(n+1)}{2}+n-1$ [Gr 1972] and  it is easy to see that there is no local  $C^\infty$-immersions  of generic smooth Riemannian $n$-manifolds  to  
$\mathbb  R^{\frac {n(n+1)}{2}-1}$. 

But for { \it no} $n\geq  2$ one can prove the existence of such immersions to $\mathbb  R^{N=\frac {n(n+1)}{2}+n-2}$  or to 
find a counterexample for $N=\frac {n(n+1)}{2}$.
  
    It is also   unclear  for which $i=2,3,....$, (if any)  all $C^\infty$-smooth Riemannian manifolds   $X=X^n$,   admit (local or global)  isometric  $C^i$-immersion  
    to $\mathbb R^{\frac {n(n+1)}{2}-1}$,   where, for all we know, this may be possible, say   for  $i\leq 0.1 \sqrt[4]{n}$,
  see [Gr 2017] for more about it.}

 The simplest non-free  $\mathcal {II}$-immersions, $X\to Y$ are those, 
 which are free on a totally geodesic hypersurface $X_0\subset X$ and also on the complement  
$X\setminus X_0$  and where  the Hessian of the second derivatives  
$\partial _{ij}f$  vanishes on $X_0$ with {\it finite order} (see 2.3.8, 3.1.8. 3.1.9 in [Gr 1986].
Then one define by induction $k$-{\it subfree} maps of $n$-manifolds $X$ to $Y$, where 

\hspace {35mm} 0-{\it subfree =free},

\hspace {-6mm} and where

\hspace {35mm}  1-{\it subfree maps} $f$

\hspace {-6mm} generalize  the above ones, namely, where    there exist 
 tangent hyperplanes $T^{n-1}_x\subset T_x(X)$   at all $x\in X$, such 
the {\it restrictions} of    $f$ to the {\it exponential  images} $X_x=\exp(T_x)\subset X$ are free  at $x$ and also on the complements $X_x$ near $x$, i.e. in the small balls: 
$$B_x(\varepsilon)\cap  (X\setminus X_x)$$ 
with    "finite order of non-freedom" on  (infinitesimal) neighbourhoods of $X_x$.

Finally,  a map   $f$   is

\hspace {35mm} {\it $k$-subfree, $k=0,1,...,n,$} 

\hspace {-6mm}  if it is  free on the above hypersurfaces $X_x\subset X $ near $x$  for all $x\in X$  and  it is 
$(k-1)$-subfree on $B_x(\varepsilon)\cap  (X\setminus X_x)$  with 
 "finite order of non-$(k-1)$-subfreedom"   on  (infinitesimal) neighbourhood of $X_x$.

{\it Clarification.} "Finite order of non-subfreedom" means that the coefficients of the  
 relevant linear differential operator $\mathscr  M=\mathscr  M_f(g)$ on
$X\setminus X_x$ are rational functions in partial  derivatives of  $f$, where these functions have their poles on $X_x$,
(see 238 in[Gr1986] and [DeL{2019]). 

This  seems to imply --I didn't truly checked this -- that 

\hspace {12mm} {\it generic immersions $X^n\to Y^N$  for $N\geq \frac {n(n+1)}{2}+n$ are $\cal II$.}

It is also plausible, that 

{\it  generic bendings of  m-subtori in the Clifford torus, as in the proof of 2,B would make
them $\cal II$ for $ \frac {m(m+1)}{2} + m(n-m)+2n-m\leq 2N.$ }

But  this, even if true only slightly improve the lower bound on $N$   in  the above  "worst dimension case", where, in fact, we expect the following.

\textbf {3.A.  Conjecture.} 
{\sf Let $X^m$ and $Y^N$  be  smooth  Riemannian manifolds, where $X$  compact and $Y$ is complete, and let  $f_0:X \to Y$ 
 be a smooth strictly short  map.
 If  
 $$N\geq  \frac {m(m+1)}{2} +1,$$
then $f_0$    admits a    $\delta$-{\it approximation}  for all  $\delta>0$  by  isometric $C^\infty$-immersions 
$f_\delta:X^m\to Y$   with
$$ curv (f(\delta(X))
 \leq \frac {1}{\delta}\left (\sqrt \frac{3m}{m+2}+\Delta(m,N)\right )+o\left (\frac {1}{\delta}\right),$$
 for the same $\Delta$ as in} 2.A, {\sf where, if 
   $$N\geq  \frac {m(m+1)}{2} +m,$$
these $f_\delta$  can be chosen $m$-{\it subfree}.}\vspace {1mm}

\hspace{45mm}{\sc Families of Maps.} \vspace {1mm}

Given a locally  defined class  $\mathscr C$ of $C^\infty $-maps $f: X\to Y$, where $Y=(Y,h)$ is a smooth Riemannian manifold,  
e.g.   $\mathscr C$ consists  of  smooth immersions, of 
$\cal II$-isometric immersions etc, the  $\delta$-approximation problem is accompanied by  a similar problem for families of maps.

 More generally, such  approximation makes sense for   maps $f$ from foliated leaf-wise Riemannian   manifolds $\mathcal X=  (\mathcal X,g) $ \footnote{This $g$ --  a smooth  leaf-wise Riemannin metric on $\mathcal X$ --  is a positive definite  differential quadratic form on the tangent bundle $\mathcal T\subset T(\mathcal X)$  to the leaves $X\subset \mathcal X$.}  to $Y$, where the restrictions of  
$ f$ to the leaves $X\subset \cal X$  are in $\mathscr C$.

\textbf {3.B. $\mathscr C$-Homotopy Approximation   Conjecture.} 
{\sf Let   $\mathcal X =(\mathcal X,g) $ be a compact manifold foliated into $m$-dimensional Riemannin leaves,
let  $Y=(Y,h)$ be complete Riemannin $N$-manifold.

Let  $\mathscr C$ be a class of smooth  leaf-wise isometric  $\cal II$-maps  and    
$\phi_0 : \mathcal X \to Y $ be  a $\mathscr C$-map.

 Let $f_0:   \mathcal X \to Y$  be a  smooth leaf-wise strictly short  $\mathscr C$-map, which is   
     homotopic to $\phi$.

 Then, at least in  in the following three cases, {\sf the map $f_0$  can be   
 $\delta$-approximated for all $\delta>0$ by smooth  leaf-wise $\mathscr C$-maps
$$f_\delta:   \mathcal X \to Y, $$
where the maps $ f_\delta$ can be joined with   $\phi_0  $ by  homotopies of  leaf-wise $\mathscr C$-maps and where
 the leaf-wise curvatures of  $f_\delta$ are bounded by

$$ curv (f(\delta(X))\leq \frac {\Xi}{\delta}+o\left (\frac {1}{\delta}\right),$$}
 %\leq \frac {\Xi}{\delta}\left (\sqrt \frac{3m}{m+2}+\Delta\right )+o\left (\frac {1}{\delta}\right),$$}
where $\Xi=\Xi(m,N, dim(\mathcal X))\leq 100$.

 {\it Case} 1. $\mathscr C$ is the class of leaf-wise {\it free} isometric 
 maps.

{\it Case} 2.  $\mathscr C$ is the class of leaf-wise $m$-{\it subfree} isometric 
 maps.}

{\it Case} 3.  $\mathscr C$ is the class of  {\it all} leaf-wise isometric 
 $ \cal II$-maps.
\vspace {1mm}

 {\it \textbf {Remarks/Questions.}} (i) {\it How big is} $\Xi$? Possibly, $\Xi$
is significantly  smaller than 100, but it is unlikely to approach  $\sqrt  3$  for $N\gg m^2$. 

Yet, it follows by the arguments in section 2  that 

{\it if  $N\gg dim(\mathcal X)^2$, then
$$\Xi \leq\sqrt { \frac{3(2dim(\mathcal X)-1) }{2dim(\mathcal X)+1}}<\sqrt 3.$$}

  Conceivably,  the correct bound on curvature needs  $\Xi\sim 3 (dim(\mathcal X))-m+1)$.} 

(ii) {\it Why Integrable?} The above  make sense for   possibly  non-integrable subbundles  
$\mathcal T\subset T(\mathcal X)$,   where the counterpart  of the  metric inducing operator sends maps $f$ to the restrictions of the forms $f^\ast(h)$ to  the subbundle $ \mathcal T\subset \mathcal X$, where  one can define the corresponding classes $\mathcal C_{\mathcal T}$  of maps $f: \mathcal X \to Y$  and where the 
$\mathcal T$-curvature of $f$ is defined as follows.

Given  a non-zero  tangent vector $\tau\in T_x(\mathcal X)$, $x\in \mathcal X$,  let
$curv_\tau(f)$ be the infimum of $Y$-curvatures at $x$  of the  
$f$-images  of the curves $C\subset \mathcal X$ which contain $x$ and are tangent to $ \tau$,
 $$curv_\tau(f)=\inf_C curv_xf(C)$$
 and
 $$curv_{\mathcal T} (f)=\sup_{\tau\in \mathcal T} curv_\tau(f).$$

Here is what one  (may be unrealistically) expects in this regard.

\textbf {3.C. Orthonormal Frame  Conjecture. } { \sf Let  $\mathcal X$ be a compact smooth Riemannian manifold and   $\Theta_i$, $i=1.,,,m\leq dim(\mathcal X)$,  be    smooth  orthonormal vector fields  on  $\mathcal X$.

 Let $Y=(Y,h)$ be a complete Riemannin $N$-manifold and 
$f_0:    \mathcal X \to Y, $ be a smooth strictly short   map.

If 
$$N\geq \frac {m(m+1)}{2}+1$$
and  the induced (pullback) vector bundle $f^\ast(T(Y))\to \mathcal X$ admits $m$ linearly independent sections,  e.g.   $f_0$ is contractible,
 then  the map $f_0$  can be   
 $\delta$-approximated for all $\delta>0$ by smooth maps $$f_\delta:   \mathcal X \to Y, $$
such the differential  images $d_f(\Theta_i)\in T(Y)$  are orthonormal with respect to $h$, 
 $$\langle d_f( \Theta_i ) , d_f(\Theta_j\rangle_h=0, \mbox { } ||d_f( \Theta_i )||_h=1,$$ and such that
$$ curv_{\mathcal T}  (f)
 \leq \frac {\Xi}{\delta}+o\left (\frac {1}{\delta}\right),$$}
for $\Xi\leq 100$,  where   $\mathcal T\subset T(\mathcal X)$ 
is  the subbundle spanned by the fields $\Theta_i$.

\vspace {1mm}

Below we state    without proof
the only confirmation  we have of 
 the conjectures 3.A and 3.B. 

\textbf  {3.D. Parametric $1D$-Approximation  Theorem.} {\sf Let $\mathcal X$ be compact  smooth Riemannin  manifold and $\mathcal T^1$ be a smooth  line  field on $\mathcal X$.
Let $Y=(Y,h)$ be a complete Riemannin manifold of dimension $N\geq 2$  and let $f_0:  \mathcal X \to Y $ be a strictly short map.
 
If  $\mathcal T^1$, regarded as a line bundle over $\mathcal X$, admits an injective homomorphism to the   induced (pullback)  bundle $f^\ast(T(Y))\to \mathcal X$, e.g.  $\mathcal T^1$ is orientable (defined by a vector field) and the map $f_0$ is contractible,     then   {\sf the map $f_0$  can be   
 $\delta$-approximated for all $ \delta>0$ by smooth 
maps 
$$f_\delta:   \mathcal X \to Y, $$ 
which are  isometric on  $\mathcal T^1$ and such that 
 the $\mathcal T^1$-curvatures of  $f_\delta$, i.e. the $Y$-curvatures  of the $f_\delta$-images of the (1-dimensional) orbits/leaves of 
$\mathcal T^1$)
, are bounded by
$$ curv _{\mathcal T^1}(f_\delta)
 \leq \frac {\Xi_1}{\delta}  +o\left (\frac {1}{\delta}\right),$$}
where $\Xi_1\leq 4.$\footnote{Probbaly, $\Xi_1< 3$   and it seems not hard to show that $\Xi_1\geq 2$. }

\textbf  {3.E. Example/Corollary.} There exists a smooth map $f: S^{2n+1} \to  
B^2(1)$, for all $n=1,2,...$,   such that  the $f$-images of all Hopf circles $S_{p}^1\subset S^{2n+1},$  $p\in\mathbb CP^n,$ 
 have equal  length $l=l_n$ $ (\leq 100n)$  and  curvatures 
 $$curv (f(S^1_p))\leq 4+\varepsilon$$
for a given $\varepsilon>0$.}

 \section {References}
 $ $ \vspace {-3mm}

[Del 2017]   R. De Leo,
{\sl Proof of a Gromov conjecture on the infinitesimal invertibility of the metric inducing operators,}
arXiv:1711.01709.
  	\vspace {2mm}

[G-J 1971] R. Greene, H. Jacobowitz,  {\sl Analytic isometric embeddings,}  Annals of Mathematics. Second Series. 93 (1): 189-204.
\vspace {2mm}

[Gr 1972]  M. Gromov, {\sl Smoothing and inversion of differential operators,}
Mathematics of the USSR-Sbornik, Volume 17, pp. 381-435.

\vspace {2mm}

[Gr 1986] M. Gromov, {\sl Partial differential relation.} 
 Springer-Verlag (1986), Ergeb. der Mat.\vspace {2mm}

[Gr 2017] M. Gromov, {\sl Geometric, Algebraic and Analytic Descendants of Nash Isometric Embedding Theorems,} Bull. Amer. Math. Soc. 54 (2017), 173-245.

\vspace {2mm}

[Gr 2022] M. Gromov,  {\sl Curvature, Kolmogorov Diameter, Hilbert
Rational Designs and Overtwisted Immersions}, arXiv:2210.13256.

\vspace {2mm}

[G-R 1970]  M.  Gromov, V. A. Rokhlin, {\sl Embeddings and immersions in Riemannian geometry}, Uspekhi Mat. Nauk, 25:5(155) (1970), 3-62.

\vspace {2mm}

[P 2023] A. Petrunin {\sl To appear.}

\end{document}